# Robust design optimization taking into account manufacturing uncertainties of a permanent magnet assisted synchronous reluctance motor

Adán Reyes Reyes
*IFP Energies nouvelles, Institut Carnot IFPEN Transports Energie*
Rueil-Malmaison, France
*& SATIE, CNRS, Paris Saclay University*
adan.reyes-reyes@ifpen.fr

Delphine Sinoquet
*IFP Energies nouvelles, Institut Carnot IFPEN Transports Energie*
Rueil-Malmaison, France
delphine.sinoquet@ifpen.fr

André Nasr
*IFP Energies nouvelles, Institut Carnot IFPEN Transports Energie*
Rueil-Malmaison, France
andre.nasr@ifpen.fr

Sami Hlioui
*SATIE, CNRS, Cergy Paris University,
Paris Saclay University*
Gif sur Yvette, France
sami.hlioui@ens-paris-saclay.fr

*Abstract*— In this paper, deterministic and robust design optimizations of a permanent magnet assisted synchronous reluctance machine were performed to increase its mean torque while reducing torque ripple. These optimizations were carried out using a surrogate model based on 2-D finite element simulations. The results of the robust optimizations, which considered manufacturing uncertainties, were compared to the deterministic optimization. The robust designs have shown not only good mean torque and torque ripple performances, but they have also shown improved robustness against design parameters uncertainties.

*Keywords— Synchronous Machines, Robust Design Optimization, Manufacturing uncertainties, Surrogate model, Finite elements analysis.*

## I. Introduction

With the increasing concerns over climate change, many measures have been adopted to reduce greenhouse gas emissions. For transportation systems, in order to replace internal combustion engine vehicles, electric and hybrid vehicles (EV, HEV) have been intensively developed. In these vehicles, the electrical machine is one of their main components.

Among the different types of electrical machines used in electric vehicles, Permanent Magnet assisted Synchronous Reluctance Machines (PMaSRMs) are one of the most used machines nowadays thanks to their good performances and their relatively low cost [1]. Unlike Surface-Mounted Permanent Magnet Synchronous Machines (SMPMSM), PMaSRMs exploit two types of torque: the hybrid torque generated using permanent magnets and, the reluctance torque which makes profit of the machine's saliency. Since SMPMSMs only generate hybrid torque, they need more permanent magnets to achieve the same torque density and tend therefore to be more expensive. However, PMaSRMs have the disadvantage of having higher torque ripple partly due to the inhomogeneous reluctance in their rotors. The design optimization phase of PMaSRMs is therefore critical to get the right balance between performances and an acceptable level of torque ripple.

Many studies have dealt with reducing torque ripple in PMaSRMs [2][3]. The optimization methodologies used in such studies can be described as deterministic since they do not consider any uncertainties on the input parameters. However, in practice, there are many discrepancies between the theoretical and real (measured) values of these parameters. These differences can be caused by manufacturing and assembly tolerances in the prototype as well as by the lack of precision on the magnetic properties of the used materials. These variabilities impact the measured performances which can diverge from those simulated in the design phase. To reduce such deviations, the parameter uncertainties should be considered in the optimization procedure.

In opposition to the deterministic optimization, robust optimization considers two types of input parameters: certain parameters also known as controllable parameters, and uncertain parameters. Controllable parameters are the same ones used in a deterministic optimization whereas uncertain parameters are specific to robust optimization techniques. This type of parameter can take varying values due to the associated uncertainties: it is then modelled by a random variable and an associated probability distribution. The presence of random input variables for the simulator leads to random output variables and then, random objective and constraint functions. Various formulations of the resulting optimization problem are proposed in the literature [4][5][6] based on expectation, probability, or quantiles of these random variables.

Reliability Based Design Optimization (RBDO) is a method used to obtain optimal and safe designs in the sense that the outputs of certain functions are inside a security domain, described by constraints. A robust or reliable design has therefore a high probability to respect these constraints. Examples of this approach can be found in [7].

Worst-case optimization considers the extreme values as objective functions and/or constraints i.e., the maximum or minimum value of the outputs at the controllable inputs caused by the uncertainty propagation [8] .

There is another very common formulation which was also adopted in this work: Robust Design Optimization (RDO). In this methodology, the expectation (average) of the objective function is optimized. To limit extreme values, a second objective based of the variances of the objective function can be added [9].

Computing these quantities requires a large sample of the uncertain input variables and thus a large number of simulations. To limit this high computational cost especially with the use of finite element simulations, meta-modelling techniques coupled with design of experiments are used to replace the costly simulations by predictions using the resulting surrogate model [10].

In this work, the results of a deterministic and two robust optimizations on a 3-phases 10-poles 60-slots PMaSRM were compared. To do so, techniques like Design Of Experiments (DOE), Finite Element Method (FEM) surrogate modelling, sensitivity analysis, quasi Monte Carlo methods and optimization algorithms were used. All the steps followed in this work will be detailed in the next sections.

## II. MACHINE TOPOLOGY AND DESIGN PARAMETERS

The machine studied in this paper is shown in Fig. 1 – a). It is a 3-phases 10-poles 60-slots PMaSRM with a Machaon rotor structure. This machine was initially designed for an EV application having a maximum torque of 430 N.m. It has an outer stator diameter of 220 mm and an active length of 200 mm. Each pole has 3 flux barriers and 7 PMs. Fig. 1 – b) shows the design parameters for the stator and for the rotor considering one layer. TABLE I lists the optimization variables with their lower and upper bounds.

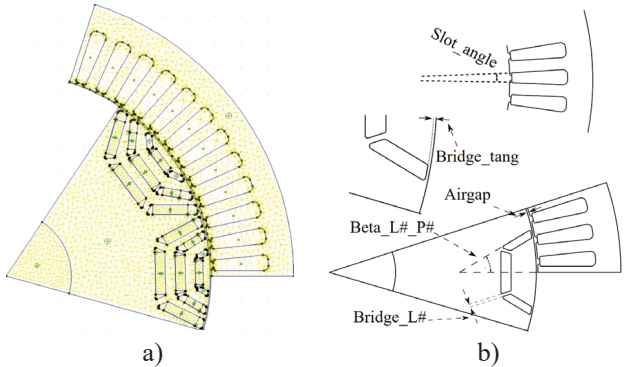

Fig. 1: a) Geometry of the PMaSRM studied in this work, b) design Parameters for one layer (# is the number of the layer).

TABLE I. Optimization variables

| Input Parameter | Description | Lower bound $x_l$ | Upper bound $x_u$ | Manufacturing Tolerance |
|---|---|---|---|---|
| Slot_angle | Stator slot width opening angle | 2.47° | 3.27° | ±0.1° |
| Beta_L1_P1 | Layer 1-Pole 1 Opening angle | 27.03° | 29.66° | ±0.33° |
| Beta_L1_P2 | Layer 1-Pole 2 Opening angle | 37.03° | 39.66° | ±0.33° |
| Beta_L2_P1 | Layer 2-Pole 1 Opening angle | 31.03° | 33.66° | ±0.33° |
| Beta_L2_P2 | Layer 2-Pole 2 Opening angle | 47.03° | 49.66° | ±0.33° |
| Beta_L3_P1 | Layer 3-Pole 1 Opening angle | 33.7° | 37° | ±0.33° |
| Beta_L3_P2 | Layer 3-Pole 2 Opening angle | 59.7° | 63° | ±0.33° |
| Airgap | Airgap width | 0.55 mm | 0.65 mm | ±0.03 mm |
| Bridge_L1 | Layer 1 radial bridge width | 2.6 mm | 2.98 mm | ±0.05 mm |
| Bridge_L2 | Layer 2 radial bridge width | 0.9 mm | 1.18 mm | ±0.05 mm |
| Bridge_L3 | Layer 3 radial bridge width | 0.5 mm | 0.62 mm | ±0.03 mm |
| Bridge_tang | Bridge between airgap and flux barriers. | 0.4 mm | 0.6 mm | ±0.05 mm |

## III. OPTIMIZATION WORKFLOW

We detail in this section the workflow used to carry out the optimizations (Fig. 2). At first, a DOE was built with the upper and lower bounds of the input parameters shown in TABLE I to fit a surrogate model for the objective functions: the mean torque and torque ripple. Secondly, these models were used to perform a global sensitivity analysis to detect the most impacting parameters on the objective functions. This will allow us to limit the number of parameters considered as uncertain. At last, and after performing the meta-model-based deterministic and robust optimizations, FEM simulations will be carried out to verify the results.

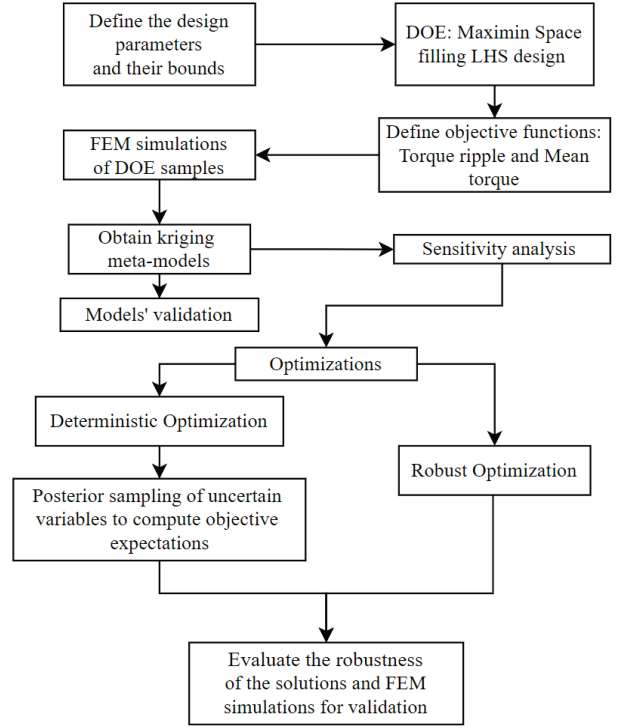

Fig. 2: Optimization workflow.

### A. Surrogate Models

To reduce computation time, surrogate models have been built for each of the objective functions. To build such models, there are three steps to follow: build a DOE, train the metamodel and check its predictivity with a test set. The chosen DOE is a maximin Latin Hypercube Sampling (LHS) as it aims to cover the search space while preserving good projection properties [11]. This DOE was built with 234 points using the bounds described in TABLE I. As for the surrogate model, it is a universal Kriging with linear trend function. For mean torque, a tensorized Matérn 5/2 covariance function has given the best predictivity. As for Torque ripple, a tensorized absolute value exponential kernel was used since the latter is not that smooth. Kriging was chosen as metamodel since it is very good at learning nonlinear objective functions and has demonstrated good performances in electrical machines optimization [12]. Finally, to evaluate the accuracy of the metamodel, a Normalized Root Mean Square Error (NRMSE) was computed on a validation test set:

$$NRSME = (||y_{real} - y_{pred}||/||y_{real}||)*100\% \quad (1)$$

where $||\cdot||$ represents the Euclidean norm. A value close to zero indicates a good model fit. The obtained *NRSME* of the kriging model for mean torque is *0.2%* and for torque ripple is 8%. These results were obtained with a train and test sets composed by 175 and 59 samples, respectively. Torque ripple depends not only on mean torque but also on torque amplitude which makes this function more difficult to model than mean torque. We consider those metamodels sufficiently accurate for performing the sensitivity analysis and the optimization procedures.

*B. Sensitivity Analysis*

To measure the impact of each input parameter on the considered outputs, a sensitivity analysis can be performed [13]. For this work, the Sobol Indices were chosen as they measure the global impact of the input variables on the output functions. The commonly used ones are the first order indices (S) and the total indices ($S_{TOTAL}$) computed with the Kriging surrogate models:

$$S_i = VAR_{X_i}(E_{X_{\sim i}}[Y|X_i])/VAR(Y) \quad (2)$$

$$S_{TOTAL,i} = E_{X_{\sim i}}[VAR_{X_i}(Y|X_{\sim i})]/VAR(Y) \quad (3)$$

where $X_{\sim i} = X_1,\cdots,X_{i-1},X_{i+1},\cdots,X_{Nx}$ and $N_x$ is the number of optimization parameters. The results of the sensitivity analysis applied to mean torque and to torque ripple are presented in Fig. 2 and Fig. 3. Only the most important inputs are displayed for better visibility. It was found that the stator slot width opening angle (Slot_angle) and the flux barrier opening angles for barriers 1 and 2 (Beta_L1_P1, Beta_L1_P2, Beta_L2_P1 and Beta_L2_P2) have the biggest impacts on the mean torque as well as on torque ripple. These 5 parameters will be then considered as uncertain variables in the robust optimization. The dispersion of these variables will be then integrated in the robust optimizations by considering a perturbation vector *U*.

*C. Deterministic and robust optimizations*

We will present in this section the results of different optimizations: a deterministic and two robust optimizations. For the deterministic optimization problem, we have:

$$min_{x \in X}[f_1(x), f_2(x)] \quad (4)$$

where $f_1$ is the opposite of the mean torque (in order to maximize it) and $f_2$ the torque ripple; *X* is the controllable parameters space defined in TABLE I. For the robust optimization problems, two formulations have been considered:

- Expectations optimization:

$$min_{x \in X(U)} [ E_U[f_1(x+U)], E_U[f_2(x+U)] ] \quad (5)$$

- Worst-case optimization:

$$min_{x \in X(U)} [ max_{u \in \Omega} f_1(x+u), max_{u \in \Omega} f_2(x+u) ] \quad (6)$$

where $X(U)= [x_{l1}-u_{u1}, x_{u1}+u_{u1}]\times\cdots\times[x_{lNx}-u_{uNx}, x_{uNx}+u_{uNx}]$

and $\Omega= [-u_{u1}, u_{u1}]\times\cdots\times[-u_{uNx}, u_{uNx}] \quad (8)$

$u_{uj}$ is the manufacturing tolerance of parameter number *j*. The tolerance for each geometrical parameter is given in

TABLE I. Based on the sensitivity study in the previous section, only 5 parameters will be considered as uncertain. Their uncertainties were considered to follow uniform distributions, i.e., $U_j \sim Unif(-u_{uj}, u_{uj})$. For parameters with no considered uncertainties, $u_u$ is simply equal to 0.

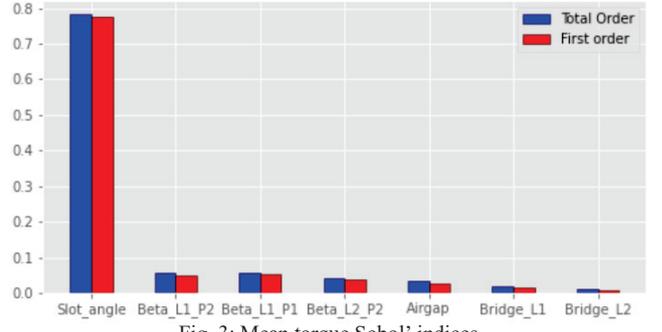

Fig. 3: Mean torque Sobol' indices

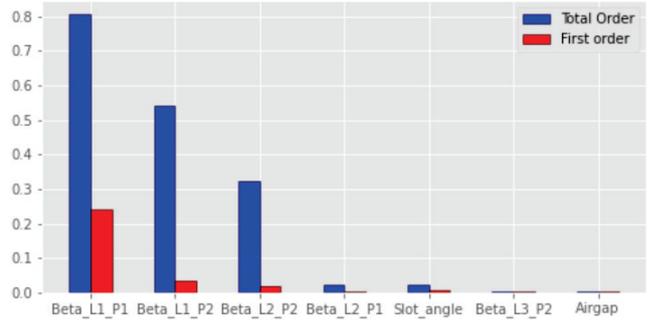

Fig. 4: Torque ripple Sobol' indices

The goal of the first robust formulation is to optimize the mean torque's and torque ripple's expectations in a Pareto sense. As for the worst-case formulation, the objective is to optimize the worst possible value of the mean torque and torque ripple caused by uncertainties. The idea is then to limit the performance degradation. To solve these optimization problems, the genetic algorithm NSGA 2 was used [14]. This algorithm has shown good performances for other studies of electrical machine optimization as in [15]. We used a DOE maximin LHS to compute samples of *x+U* to calculate the objective functions' expectations with a quasi-Monte Carlo method. When it comes to the Worst-case formulation, we have two options: computing samples of *x+U* and take the maximum value of these samples as an estimator of $max_{u \in \Omega} f(x+u)$ *or* obtaining the absolute maximum value with an optimization algorithm. In this work, the latter was adopted using a Particle Swarm Optimization (PSO) [16] algorithm for the embedded mono-objective optimizations.

Fig. 5 shows a comparison between the deterministic (blue) and the Expectations optimization (red) Pareto fronts. The expected performances of the deterministic Pareto front have been reevaluated (pink): The design variables were perturbed by adding sampled values of the uncertain variables. These expected values represent the average mean torque and average torque ripple for each machine obtained by the deterministic Pareto optimization considering a posteriori uncertainty on the input parameters. As we can notice, an optimal deterministic design does not guarantee its performances when manufacturing tolerances are considered: even though the front of the deterministic optimization (blue) shows better performances than the Expectations optimization one (red) at first glance, the torque ripple expectation values (pink) are larger than the robust front indicating that the proposed deterministic solution designs are more sensitive to uncertainties.

The pareto front of the Worst-case optimization (green) is presented in red in Fig. 6. As in Fig. 5, the deterministic pareto front was also added (blue). The worst-case performances of the deterministic Pareto front have been evaluated (light green) thanks to a posteriori uncertainty on the input parameters and PSO maximization. Once again, these results show the importance of a robust optimization in limiting the performance degradation that a sample of machines can have.

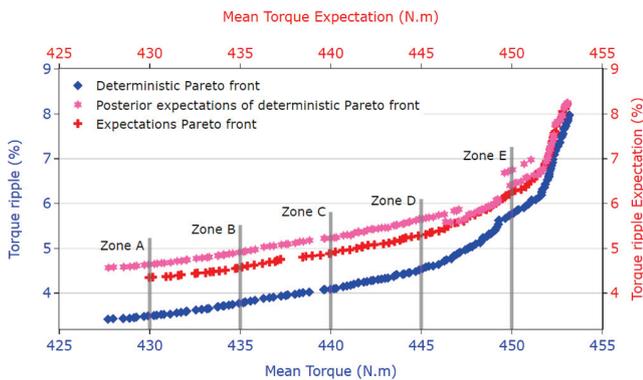

Fig. 5: Pareto fronts obtained by deterministic optimization (blue: Mean Torque and Torque ripple) and robust optimization (red: Mean Torque expectation and Torque ripple expectation). Posterior perturbations of solutions of deterministic optimization (pink: Mean Torque expectation and Torque ripple expectation). Dark gray zones highlight points with similar Mean Torque expectation values (430 ± 0.1, 435 ± 0.1, 440 ± 0.1, 445 ± 0.1 and 450 ± 0.1 N.m).

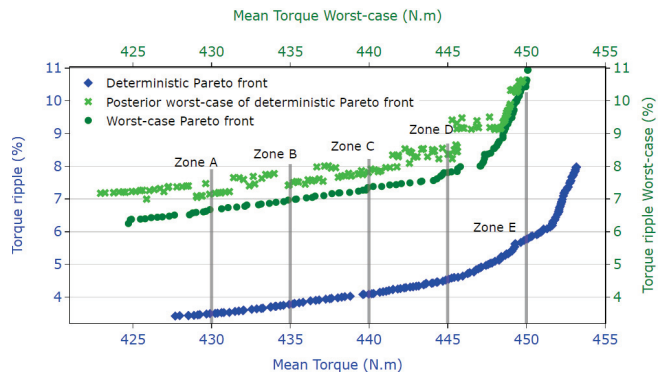

Fig. 6: Pareto fronts obtained by deterministic optimization (blue: Mean Torque and Torque ripple) and robust optimization (green: Mean Torque worst-case and Torque ripple worst-case). Posterior perturbations of solutions of deterministic optimization (light green: Mean Torque worst-case and Torque ripple worst-case). Dark gray zones highlight points with similar worst case Mean Torque values (430 ± 0.1, 435 ± 0.1, 440 ± 0.1, 445 ± 0.1 and 450 ± 0.1 N.m).

To go deeper into this analysis, we empirically compared the distribution of different designs. For this purpose, we show in Fig. 7 and Fig. 8 boxplots of a subset of points selected from the deterministic and robust Pareto fronts shown in Fig. 5 and Fig. 6, respectively. Each pair of boxplots represents a comparison of the distribution of torque ripple values between a determinist machine (blue) and a robust machine (red, green) falling in one of the zones (A, B, C, D and E). These boxplots show the values of q1 (the value for which there is 25% of the samples below it), q2, which is the median of the sample (there is 50% of the sample below it) and q3 (there is 75% of the samples below it) as represented in Fig. 7.

We can observe in this figure that for each pair of machines, the robust one shows better overall performance than the deterministic one. For example, for the machines in zone A, we can notice that the median of the deterministic machine is almost the same as the q3 value of the robust machine (4.5%). This means that there is a 50% chance for the deterministic design and 75% chance for the robust design of having a torque ripple value smaller than 4.5%. Another remarkable result is that the expectation and Standard Deviation (STD) of the torque ripple associated with the robust optimization solution (4.3%, 0.35%, respectively) outperform the expectation and STD obtained by the posterior analysis of the deterministic solution (4.6%, 0.56%, respectively). These results stress the importance of robust optimizations when dealing with uncertainties.

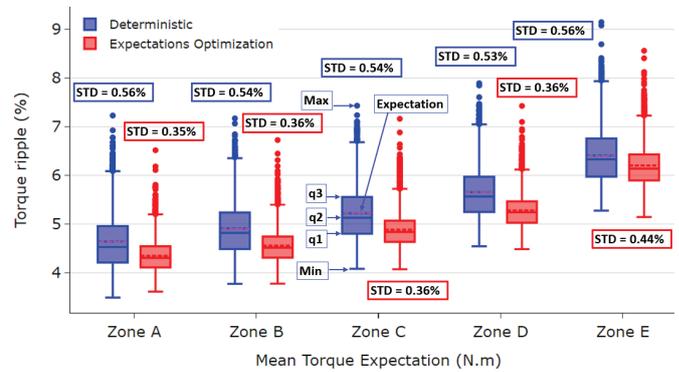

Fig. 7: Boxplots showing comparisons of predicted torque ripple values between deterministic (blue) and robust (red) machines from Fig. 5 with similar predicted Mean Torque Expectation values.

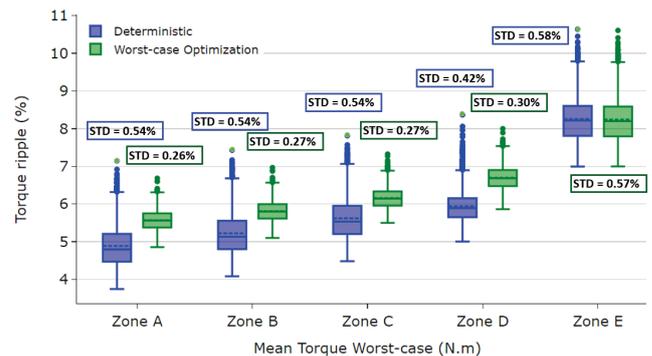

Fig. 8: Boxplots showing comparisons of predicted torque ripple values between deterministic (blue) and robust (green) machines from Fig. 6 with similar predicted Mean Torque Worst-case values.

Fig. 8 also shows pairs of boxplots comparing predicted torque ripple values between deterministic (blue) and robust (green) machines from Fig. 6 with similar predicted Mean torque worst-case values. For all the selected zones, the worst-case torque ripple value of a robust design is lower (or equal) than a deterministic one. Besides this, robust solutions also have lower STD values than deterministic solutions, especially for zones with low mean torque worst case. In zone A for example, the robust machine has a STD of 0.26% and a worst-case torque ripple of 6.7% compared to 0.54% and 7.2%, respectively for its deterministic counterpart. Both optimizations lead practically to the same design for high values of mean torque worst case like in zone E. This can also be seen in Fig. 5 and Fig. 6, with the Pareto fronts of the deterministic and robust optimizations getting very close with increasing torque.

Although the worst-case optimization has allowed to limit the performance degradation of the least performant machine in a sample, it has led to worse expectance values on torque ripple. Instead of using the worst-case as an objective function, it could be used as a constraint in a constrained optimization problem while still using the expectations as an objective. Such formulation allows to have good machines samples while limiting the worst performances we can have. The main disadvantage of the robust machines through a worst-case optimization is that there is still a big probability of having a torque ripple value greater than the deterministic machine (min, q1, q2, q3). For instance, the deterministic design in zone A has around 75% chance of obtaining a torque ripple lower than 5.2%, while the q1 value of the robust design is 5.4%.

In order to compare these two robust formulations, we propose in Fig. 9 new boxplots comparing designs issued from both optimizations and belonging to the same zones. This time, a zone corresponds to a mean torque expectation value. The robust optimization based on expectations provides not only better results in terms of expectations, but also in terms of q1, q2 and q3. For example, in zone A, a red machine has a 75% chance of having a torque ripple value outperforming any of the machines in the green sample. This is because the minimum value (4.6%) of the sample almost coincides with the q3 value of the blue machine (4.5%). Comparisons of worst-case values between boxplots from these two optimizations are inconclusive (no significant advantage in terms of worst-case torque ripple for the worst-case optimization). It can also be noticed that the worst-case optimization leads to lower STD torque ripple values. This can be explained by two reasons. First, in a torque ripple worst-case optimization, we are interested in limiting the extreme upper values which tends to bring closer the performances of the boxplots' samples. Second, in such optimization, reducing extreme lower values is not contemplated. Due to that, we risk obtaining non-optimized extreme lower values. This has been the case in our optimizations where the minimum values of the green samples are higher than the red ones. The combination of these 2 facts leads to reduced STD values.

*D. Results verification*

All the results presented in the previous sections were based on the surrogate models. We will therefore verify in this section some results using FEM simulations. Fig. 10 and Fig. 11 show the same boxplots in Fig. 7, Fig. 8, however, this time, FEM simulations were used to determine the performances. Only zones A and E are represented in these figures. The STD values were added on these figures.

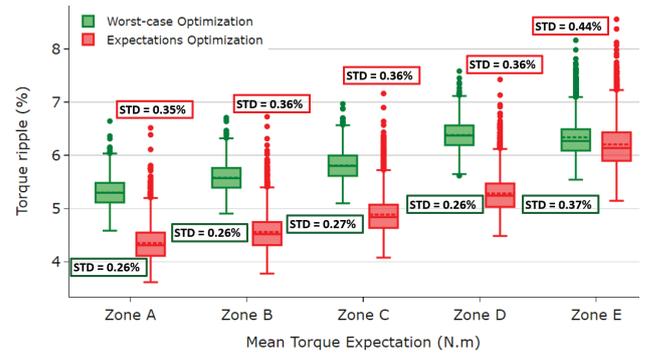

Fig. 9: Boxplots showing comparisons of predicted torque ripple values between machines from worst-case (green) and expectations (red) optimizations with similar predicted Mean Torque Expectation values.

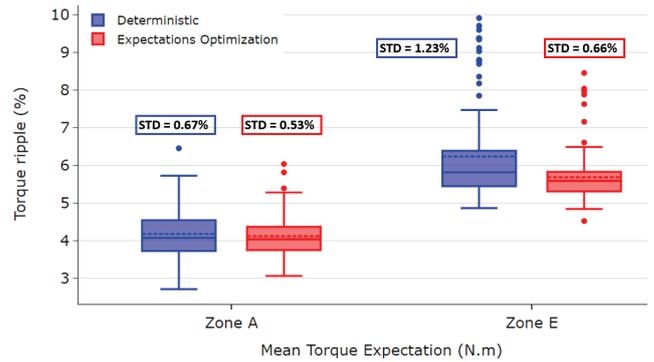

Fig. 10: Boxplots showing comparisons of FEM simulations' torque ripple values between deterministic (blue) and robust (red) machines with similar predicted Mean Torque/ Mean Torque Expectation values.

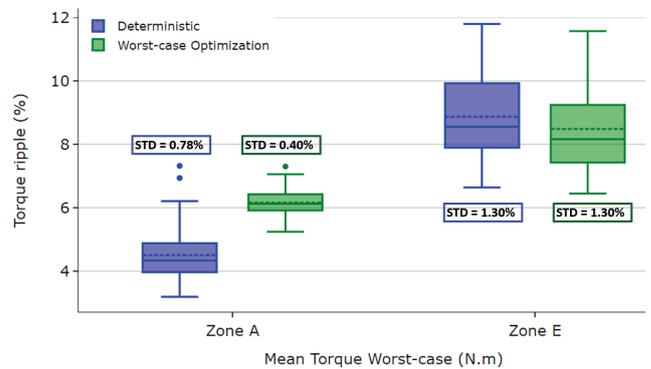

Fig. 11: Boxplots showing comparisons of FEM simulations' torque ripple values between deterministic (blue) and robust (green) machines with similar predicted Mean Torque Worst-Case values.

The results from Fig. 10 confirm what we have already seen in Fig. 7 using the meta-model. For both zones A and E, the robust optimization presents similar or better solutions than the deterministic optimization in terms of robustness. For zone A, both machines have practically the same mean torque ripple with a smaller STD for the robust design. As for zone E, the robust design presents a particularly lower maximum value, which means that in a worst-case scenario, the robust design will have a torque ripple of 8.4 % vs. 9.9 % for the deterministic design. The q3 value of the robust design is also equal to the median of the deterministic one (5.8 %), meaning that 75 % of the produced machines using the robust design would have a torque ripple lower than 5.8 % vs. only 50 % for the deterministic design. As also seen in Fig. 8, using the optimization formulation with the worst-case scenario tends to reduce the STD of the robust designs compared to the deterministic ones, especially at low torque. This is also confirmed by FEM simulations in Fig. 11. For the machines in zone A, the robust design has a notably lower STD compared to the deterministic one (0.4 % vs 0.78 %). As for the worst-case torque ripple value, both designs show similar performances.

Regarding the precision of the meta-models compared to FEM simulations, some differences can be noticed. For example, the NRMSE between predicted values in Fig. 7 and FEM simulations in Fig. 10 is 15.1% for torque ripple and 0.28% for mean torque. This lack of precision, especially for torque ripple, is somehow expected with the strategy used to create fixed meta-models and the difficulty to fit a torque ripple meta-model as seen in section III.A. A vast number of simulations is needed in this case to obtain acceptable level of accuracy. An alternative approach could be to use an adaptive strategy to update the surrogate model with additional simulations during the optimization. By using this approach, additional FEM simulations are performed only in promising points in the search space, limiting computational time and increasing the precision for optimal designs. Finally, we note that surrogate-based optimizations require much less computational time compared to FEM-based optimizations. For example, 17 hours were needed to finalize the DOE simulations *(Intel(R) Xeon(R) W-2195 CPU @ 2.30 GHz and 18 cores)*. The deterministic optimization using the meta-model took only 13 seconds for 300 iterations and 150 particles. Doing the same optimization with only FEM simulations would have taken around 3 months to complete.

IV. CONCLUSIONS

We have presented in this paper a comparison between three different optimizations performed on a PMaSRM to minimize its torque ripple and maximize its mean torque. The first optimization used a deterministic formulation and the other ones used robust formulations which take input parameters' uncertainties into account by setting the objective functions as expectations and worst-case values. In order to reduce computation time, surrogate models have been built for each of the objective functions. These surrogate models have been also used to perform a sensitivity analysis to detect the most impacting input parameters. Objective functions' expectations were computed with a quasi-Monte Carlo scheme while worst-cases where calculated with a PSO algorithm. It should be noted that while the predicted values of mean torque were consistent with FEM simulations, some differences were observed for torque ripple. This problem will be addressed in future projects by using an adaptive strategy to update the surrogate models with additionnal simulations during the optimization. Nevertheless, the FEM simulations have confirmed the tendency of the predicted results using the metamodels. The comparison of Pareto fronts has shown that robust solutions outperform deterministic solutions in terms of robustness criteria. The expectations as well as the STD values using robust optimizations were better that for the deterministic optimization. This shows the importance of developing new techniques of optimization of electrical machines when dealing with manufacturing uncertainties.